\documentclass[11pt]{article}

\usepackage{times,amsfonts}
\usepackage[a4paper]{geometry}
\newtheorem{thm}{Theorem}[section]
\newtheorem{lem}[thm]{Lemma}
\newenvironment{proof}{\noindent {\bf Proof}.}{\hfill$\Box$\\[-5mm]}

\newcommand{\NN}{\mathbb{N}}
\newcommand{\RR}{\mathbb{R}}
\newcommand{\DD}{\mathbb{D}}

\newcommand{\EE}{\mathbf{E}}
\newcommand{\PP}{\mathbf{P}}

\newcommand{\One}{{\mathrm{1} \kern -0.27em \mathrm{I}}}

\newcommand{\defeg}{\mathrel :=}

\newcommand{\ut}{{\tilde u}}
\newcommand{\lin}{\mathbf{K}}
\newcommand{\lint}{\mathcal{K}}

\begin{document}

\title{\vskip-10mm
Out of equilibrium functional central limit theorems
for a large network where customers join the shortest of several queues}

\author{\sc
Carl Graham~\footnote{CMAP, 
{\'E}cole Polytechnique, 91128 Palaiseau, France. UMR CNRS 7641.
{\tt carl@cmapx.polytechnique.fr} 
December 1, 2003.}
}

\date{}

\maketitle

\begin{quote}\small  \vskip -8mm
{\em Abstract\/}. Customers arrive at rate $N\alpha$
on a network of $N$ single server infinite buffer queues,
choose $L$ queues uniformly, join the shortest one, and are served
there in turn at rate $\beta$.
We let $N$ go to infinity.
We prove a functional central limit theorem (CLT) for the tails 
of the empirical measures of the queue occupations,
in a Hilbert space with the weak topology, with limit
given by an Ornstein-Uhlenbeck process. 
The {\em a priori\/} assumption is that the initial data converge.
This completes a recent functional CLT in equilibrium 
in Graham~\cite{Graham:03} for which convergence for the initial data 
was not known {\em a priori\/}, but was deduced {\em a posteriori\/} from the 
functional CLT.


\vskip -25mm
\end{quote}

\section{Introduction}

\subsection{The queuing model}

We continue the asymptotic study 
for large $N$ and fixed $L$
initiated in
Vvedenskaya et al.~\cite{Vved:96}
of a Markovian network constituted of $N$ single server infinite buffer  queues. 
Customers arrive at rate $N\alpha$, are allocated
$L$ distinct queues uniformly at random, and join the shortest, ties
being resolved uniformly at random. Service is at rate $\beta$.
Arrivals, allocations and services 
are independent.
The interaction structure depends on sampling 
from the empirical measure of $L$-tuples of queue
states; in statistical mechanics terminology, this constitutes
$L$-body mean-field interaction. 

Let $X^N_i(t)$ be the length of queue $i$ at time $t\ge0$.
The process $(X^N_i)_{1\leq i\leq N }$ is Markov, 
its empirical measure $\mu^N  ={1 \over N}\sum_{i=1}^N \delta_{X^N_i}$ 
has samples in $\mathcal{P}(\DD(\RR_+,\NN) )$, and its marginal process 
$(\mu^N_t)_{t\ge0}$ has sample paths  
in $\DD\!\left(\RR_+,\mathcal{P}(\NN)\right)$. We are interested 
in the tails of the marginals, and consider
\[
\mathcal{V} 
= \left\{ 
(v(k))_{k\in\NN} : v(0)=1,\ v(k)\geq v(k+1),\ \lim v =0 
\right\}\,,
\qquad 
\mathcal{V}^N = \mathcal{V} \cap {1\over N}  \NN^\NN\,,
\]
with the uniform topology (which coincides here with the product topology)
and  the process $R^N = (R^N_t)_{t\ge0}$ with 
sample paths in $\DD\!\left(\RR_+,\mathcal{V}^N\right)$
given by
\[
R^N_t(k) = {1 \over N}\sum_{i=1}^N \One_{X^N_i(t)\ge k}\,.
\] 

The processes $(\mu^N_t)_{t\ge0}$ and $(R^N_t)_{t\ge0}$ 
are in relation through
$p \in \mathcal{P}(\NN) \longleftrightarrow v \in \mathcal{V}$
for $v(k) = p[k,\infty )$ and $p\{k\} = v(k) - v(k+1)$ for
$k$ in $\NN$. This classical homeomorphism 
maps the subspace of probability measures with 
finite first moment
onto  $\mathcal{V}\cap \ell_1$,  corresponding to
a finite number of customers in the network.
The symmetry structure implies that these processes 
are Markov.

\subsection{Laws of large numbers}

Let $c_0^0$ and $\ell_p^0$ for $p\ge1$ denote
the subspaces  of sequences vanishing at $0$
of the classical sequence spaces $c_0$ (with limit 0) and
$\ell_p$ (with summable absolute $p$-th power). 
We define mappings with values in $c^0_0$ given
for $v$ in $c_0$ by
\begin{equation}
\label{F}
F_+(v)(k) = \alpha\!\left( v(k-1)^L -v(k)^L  \right),
\quad 
F_-(v) (k) = \beta  (v(k) - v(k+1))\,,
\qquad 
k\ge1\,,
\end{equation}
and $F = F_+ - F_-$, and the nonlinear differential equation 
$\dot u = F(u)$
on $\mathcal{V}$, explicitly for $t\ge0$ 
\begin{equation}
\label{is}
\dot u_t(k)
= F(u_t)(k)
= \alpha\! \left( u_t(k-1)^L-u_t(k)^L \right)
- \beta\! \left( u_t(k) - u_t(k+1)\right)\,,
\qquad k\geq 1\,.
\end{equation}

This corresponds to (1.6) in Vvedenskaya et al.~\cite{Vved:96} 
(with arrival rate $\lambda$ and
service rate $1$) and  (3.9) in  Graham~\cite{Graham:00} 
(with arrival rate $\nu$ and service rate $\lambda$).
Theorem~1~(a) in \cite{Vved:96}
and Theorem~3.3 in \cite{Graham:00} yield that
there  exists a unique solution $u=(u_t)_{t\ge0}$
taking values in $\mathcal{V}$ for (\ref{is}), 
which is continuous, 
and if $u_0$ is in $\mathcal{V}\cap\ell_1$
then $u$ takes values in $\mathcal{V}\cap\ell_1$.

A functional law of large numbers (LLN) 
for converging initial data follows from
Theorem~2 in \cite{Vved:96}. We give below a result contained 
in Theorem~3.4 in \cite{Graham:00}.

\begin{thm}
\label{lln}
Let $(R^N_0)_{N\ge L}$ converge in probability  to
$u_0$ in $\mathcal{V}$.
Then  $(R^N)_{N\ge L}$ converges in probability in $\DD(\RR_+,\mathcal{V})$
to the unique solution $u=(u_t)_{t\ge0}$ 
starting at $u_0$ for (\ref{is}).
\end{thm}

The networks are stable for
$\rho = \alpha / \beta <1$. Then
Theorem 1 (b) in Vvedenskaya et al.~\cite{Vved:96}  yields
that (\ref{is}) has a globally stable point $\ut$ in $\mathcal{V}\cap\ell_1$
given by $\ut(k)= \rho^{(L^k -1) / (L-1)}$. 
A functional LLN in equilibrium
for $(R^N)_{N\ge L}$
with limit $\ut$ follows by a
compactness-uniqueness method validating the inversion of limits 
for large sizes and large times, 
see Theorem 5 in \cite{Vved:96} 
and Theorem 4.4 in \cite{Graham:00}.

The results of \cite{Vved:96} are extended in Graham~\cite{Graham:00},
in particular to LLNs and
propagation of chaos results on path space. 
Theorem~3.5 in \cite{Graham:00}  gives convergence bounds in variation norm 
for the chaoticity result on $[0,T]$ for  $(X^N_i)_{1\leq i \leq N}$ for
 $(X^N_i(0))_{1\leq i \leq N}$   i.i.d. of law $q$, 
using results in Graham and M{\'e}l{\'e}ard~\cite{Meleard:94}. These bounds can be somewhat
extended  for initial data satisfying appropriate a priori controls, but 
behave exponentially badly for large $T$. 
 
\subsection{Central limit theorems}

Graham~\cite{Graham:03} and the present paper seek asymptotically tight
rates of convergence and confidence intervals, and study the fluctuations
around the LLN limits. 
For $R^N_0$ in  $\mathcal{V}^N$ and $u_0$ in $\mathcal{V}$
we consider the process $R^N$, the solution $u$ for (\ref{is}),
and the  process $Z^N = (Z^N_t)_{t\ge0}$
with values in $c^0_0$ given by
\begin{equation}
\label{flu}
 Z^N = \sqrt{N} (R^N - u)\,.
\end{equation}

Graham~\cite{Graham:03} focuses on the {\em stationary} regime for $\alpha<\beta$
defining the initial data {\em implicitly\/}:
the law of $R^N_0$ is the invariant law for $R^N$
and $u_0=\ut$. 
The main result in \cite{Graham:03} is Theorem~2.12, 
a functional central limit theorem (CLT) 
in equilibrium for
$(Z^N)_{N\ge L}$ with limit 
a stationary Ornstein-Uhlenbeck process.
This {\em implies\/} a CLT under the
invariant laws for $(Z^N_0)_{N\ge L}$  
with limit the invariant law for this Gaussian process, an 
important result which seems very difficult to obtain 
directly. 
The proofs actually
involve appropriate {\em transient\/} regimes, ergodicity, 
and fine studies of the long-time behaviors, in particular
a global exponential stability result for the nonlinear dynamical system (\ref{is}) 
using intricate comparisons with
linear equations and their spectral theory.

We complete here the study in \cite{Graham:03} and derive a
functional CLT in relation to Theorem~\ref{lln},
for the Skorokhod topology on
Hilbert spaces with the weak topology, for a wide class of 
$R^N_0$ and $u_0$ under the {\em assumption\/} that
$(Z^N_0)_{N\ge L}$ converges in law (for instance satisfies a CLT). 
This covers without constraints on $\alpha$ and $\beta$ many transient regimes
with {\em explicit\/} initial conditions, such as i.i.d.\ queues 
with common law appropriately converging as $N$ grows. 
Section \ref{Sfclt} introduces in turn the main notions and results, 
and Section \ref{Sproofs} leads progressively to the  proof of the functional CLT
by a compactness-uniqueness method.

\section{The functional central limit theorem}
\setcounter{equation}{0}
\label{Sfclt}

For a sequence $w=(w(k))_{k\ge 1}$ such that $w(k)>0$
we define the Hilbert spaces
\[
L_2(w) = \biggl\{ x \in \RR^\NN : x(0)=0\,,\ 
\Vert x \Vert_{L_2(w)}^2 = \sum_{k\ge1} \biggl({x(k)\over w(k)}\biggr)^2 w(k)
= \sum_{k\ge1} x(k)^2 w(k)^{-1} <\infty
\biggr\}
\]
of which the elements are considered as measures identified with 
their densities with respect to the reference measure $w$. 
Then $L_1(w) = \ell_1^0$ and if
$w$ is summable then
$\Vert x \Vert_1 \le \Vert w \Vert_1^{1/2}\Vert x \Vert_{L_2(w)}$
and $L_2(w) \subset \ell_1^0$.
For bounded $w$ we have the Gelfand triplet
$
L_2(w) \subset \ell_2^0  \subset L_2(w)^* = L_2(w^{-1})
$.

Also, $L_2(w)$ is an $\ell_2$ space with weights, and 
we consider the $\ell_1$ space with same weights
\[
\ell_1(w) = \biggl\{ x \in \RR^\NN : x(0)=0\,,\ 
\Vert x \Vert_{\ell_1(w)} = \sum_{k\ge1} |x(k)| w(k)^{-1} <\infty
\biggr\}\,.
\]
Clearly $x \in L_2(w) \Leftrightarrow x^2 \in \ell_1(w)$.
The operator norm of the
inclusion  $\mathcal{V} \cap \ell_1(w) \hookrightarrow \mathcal{V} \cap L_2(w)$
is
bounded by 1 since
$\Vert x \Vert_{L_2(w)}^2 = \Vert x^2 \Vert_{\ell_1(w)}
\le \Vert x \Vert_{\ell_1(w)}$ for $\Vert x \Vert_\infty \le 1$.

In the sequel we assume that $w=(w_k)_{k\ge1}$ satisfies the condition that 
\begin{equation}
\label{compw}
\exists\, c,d>0 
\,:\;
c  w(k+1)\le  w(k) \le d w(k+1) \;\mbox{for}\; k\ge1\,.
\end{equation}
This holds for $\theta>0$  for the geometric sequence $(\theta^k)_{k\ge1}$, 
yielding quite  strong norms for $\theta<1$.

\begin{thm}
\label{isw}
Let $w$ satisfy (\ref{compw}).
Then in $\mathcal{V}$ the mappings
$F$, $F_+$ and $F_-$ are Lipschitz for 
the $L_2(w)$ and the $\ell_1(w)$ norms. Existence and uniqueness holds
for (\ref{is}) in  $\mathcal{V} \cap L_2(w)$ and in $\mathcal{V} \cap \ell_1(w)$.
\end{thm}

\begin{proof}
We give the proof for $\ell_1(w)$, the proof for $L_2(w)$ being similar 
(see Theorem 2.2 in Graham~\cite{Graham:03}).
The identity $x^L-y^L = (x-y)(x^{L-1} + x^{L-2}y +\cdots + y^{L-1})$ yields
\begin{eqnarray*}
\left|u(k-1)^L - v(k-1)^L\right| w(k)^{-1}
&\le& 
\left|u(k-1) - v(k-1)\right| L d w(k-1)^{-1}\,,
\\
\left|u(k)^L - v(k)^L\right| w(k)^{-1}
&\le& 
\left|u(k) - v(k)\right| L  w(k)^{-1}\,,
\\
\left|u(k+1) - v(k+1)\right|  w(k)^{-1}
&\le&
\left|u(k+1) - v(k+1)\right|  c^{-1}w(k+1)^{-1}\,,
\end{eqnarray*}
hence 
$
\Vert F_+(u) - F_+(v) \Vert_{\ell_1(w)} 
\le \alpha L (d+1)
\Vert u - v  \Vert_{\ell_1(w)}
$ and
$
\Vert F_-(u) - F_-(v) \Vert_{\ell_1(w)} 
\le \beta  (c^{-1}+1)
\Vert u - v  \Vert_{\ell_1(w)}
$. 
Existence and uniqueness follows using a Cauchy-Lipschitz method.
\end{proof}

The linearization of (\ref{is}) around a particular solution $u$ in $\mathcal{V}$
is the linearization of the equation satisfied by $z=g -u$ 
where $g$ is a generic solution for (\ref{is}) in $\mathcal{V}$, 
and is given for $t\ge0$ by
\begin{equation}
\label{linis}
\dot z_t = \lin(u_t) z_t
\end{equation} 
where for $v$ in $\mathcal{V}$ the linear operator 
$\lin(v) : x \mapsto \lin(v) x$ on $c_0^0$ is given by
\begin{equation}
\label{defK}
\lin(v) x (k) =
\alpha L v(k-1)^{L-1} x(k-1)
-(\alpha L v(k)^{L-1} + \beta) x(k)
+  \beta x(k+1)\,,\quad 
k\ge1\,.
\end{equation}
The infinite matrix in the canonical basis 
$(0,1,0,0 \ldots)$, $(0,0,1,0 \ldots)$, \dots\, is given by
\\
\[
\pmatrix{
 -\left(\alpha L v(1)^{L-1} + \beta\right) &  \beta & 0  &\cdots 
\cr
 \alpha L v(1)^{L-1} & -\left( \alpha L v(2)^{L-1} + \beta\right) &  \beta   & \cdots
\cr
 0 &  \alpha L v(2)^{L-1} & -\left( \alpha L v(3)^{L-1} + \beta\right)  & \cdots
\cr
 0 &  0 & \alpha L v(3)^{L-1} & \cdots
\cr
 \vdots & \vdots & \vdots  & \ddots 
}
\]
\\
and $\lin(v)$ is the adjoint of the
the infinitesimal generator of a sub-Markovian birth and death process. The
spectral representation of Karlin and McGregor~\cite{Karlin:57a} was a key tool 
in Graham~\cite{Graham:03}, but here it varies in time and introduces no true 
simplification.

Let  $(M(k))_{k\in\NN}$ be independent
real continuous centered Gaussian martingales, determined in law by their
deterministic Doob-Meyer brackets given for $t\ge 0$ by
\begin{equation}
\label{limdoob}
 \langle M(k) \rangle_t = \int_0^t
\left\{ F^{\vphantom{N}}_+(u_s)(k) + F^{\vphantom{N}}_-(u_s)(k) \right\}ds\,. 
\end{equation}
The processes $M = (M(k))_{k\ge0}$ and
$\langle M \rangle = \left(\langle M(k) \rangle\right)_{k\in\NN}$
have values in $c^0_0$.

\begin{thm}
Let $w$ satisfy (\ref{compw}) and $u_0$ be in $\mathcal{V} \cap \ell_1(w)$.
Then the Gaussian martingale $M$ is square-integrable in $L_2(w)$.
\end{thm}

\begin{proof}
We have $\EE\!\left(\vphantom{M^N}\smash{\Vert M_t \Vert_{L_2(w)}^2}\right)
= \EE\!\left(\vphantom{M^N}\smash{\Vert \langle M\rangle_t \Vert_{\ell_1(w)}}\right)$ 
and we conclude using (\ref{limdoob}), Theorem \ref{isw},
and uniform bounds in 
$\ell_1(w)$ on $(u_s)_{0\le s \le t}$  in function of $u_0$
given by the  Gronwall Lemma.
\end{proof}

The limit Ornstein-Uhlenbeck equation for the fluctuations is 
the inhomogeneous affine SDE given for $t\ge0$ by
\begin{equation}
\label{ousde}
Z_t = Z_0 + \int_0^t \lin(u_s) Z_s \,ds + M_t
\end{equation}
which is a perturbation of (\ref{linis}).
A well-defined solution is called an Ornstein-Uhlenbeck
process. 

In equilibrium $u=\tilde u$ and
setting $\lint = \lin(\ut)$ and using (\ref{F}) and $F_+(\ut)=F_-(\ut)$
yields the simpler and more explicit formulation in Section~2.2 in 
Graham~\cite{Graham:03}. 
We recall that strong (or pathwise) uniqueness implies weak uniqueness, and
that $\ell_1(w) \subset L_2(w)$.

\begin{thm}
\label{Kbdd}
Let the sequence  $w$ satisfy (\ref{compw}). 
\\
(a) 
For $v$ in $\mathcal{V}$, the operator $\lin(v) $ is bounded in $L_2(w)$,
and its operator norm is uniformly bounded.
\\
(b) Let $u_o$ be in $\mathcal{V} \cap L_2(w)$. Then in $L_2(w)$ there
is a unique solution 
$z_t = \mathrm{e}^{\int_0^t \lin(u_s)\,ds} z_0$ for (\ref{linis}) and
strong uniqueness of solutions holds for (\ref{ousde}).
\\
(c) Let $u_o$ be in $\mathcal{V} \cap\ell_1(w)$.
Then in $L_2(w)$ there
is a unique strong solution 
$
Z_t = \mathrm{e}^{\int_0^t \lin(u_s)\,ds} Z_0 
+ \int_0^t \mathrm{e}^{\int_s^t \lin(u_r)\,dr} dM_s
$
for (\ref{ousde})
and if $\EE\Bigl(\Vert Z_0 \Vert_{L_2(w)}^2\Bigr) < \infty$ then 
$\EE\Bigl(\sup_{t\le T}\Vert Z_t \Vert_{L_2(w)}^2\Bigr) < \infty$.
\end{thm}

\begin{proof}
Considering (\ref{defK}), $v\le1$, convexity bounds, and (\ref{compw}),  we have 
\begin{eqnarray*}
\Vert \lin(v) x \Vert_{L_2(w)}^2
&\le&
2(\alpha L + \beta)
\sum_{k\ge1}
\left(
\alpha L x(k-1)^2 d w(k-1)^{-1} + 
(\alpha L + \beta) x(k)^2 w(k)^{-1}
\right.
\\
&&\kern 30mm\left. 
{} + \beta x(k+1)^2 c^{-1}  w(k+1)^{-1}
\right)
\\
&\le&
2(\alpha L + \beta)(\alpha L (d+1) + \beta (c^{-1} +1))
\Vert x \Vert_{L_2(w)}^2
\end{eqnarray*}
and (a) and (b) follow,
the Gronwall Lemma yielding uniqueness.
Under the assumption on $u_0$ in (c) the martingale $M$ is square-integrable
in $L_2(w)$. If $\EE\Bigl(\Vert Z_0 \Vert_{L_2(w)}^2\Bigr) < \infty$ then 
the formula for $Z$ is well-defined, solves the SDE,
and the Gronwall Lemma yields 
$\EE\Bigl(\sup_{t\le T}\Vert Z_t \Vert_{L_2(w)}^2\Bigr) < \infty$, else
for any $\varepsilon>0$
we can find $r_\varepsilon<\infty$ such that 
$\PP\!\left(\Vert Z_0 \Vert_{L_2(w)} \le r_\varepsilon\right) > 1- \varepsilon$,
and a localization procedure using pathwise uniqueness yields existence.
\end{proof}

Our main result is 
the following functional CLT.
We refer to Jakubowski~\cite{Jakubowski:86} for the Skorokhod topology
for non-metrizable topologies.
For the weak topology of a 
reflexive Banach space, the relatively compact sets are the 
bounded sets for the norm, see Rudin~\cite{Rudin:73} Theorems 1.15~(b),
3.18, and 4.3. Hence,  $B(r)$ denoting the  closed ball
centered at 0 of radius $r$, 
a set $\mathcal{T}$ of probability measures
is tight if and only if for all $\varepsilon>0$ there exists 
$r_\varepsilon < \infty$ such that $p(B(r_\varepsilon))> 1 - \varepsilon$
uniformly for $p$ in $\mathcal{T}$, which is the case if $\mathcal{T}$
is finite.

\begin{thm}
\label{fclt}
Let $w$ satisfy (\ref{compw}). Consider 
$L_2(w)$
with its weak topology and
$\DD(\RR_+, L_2(w))$ with the corresponding Skorokhod
topology.
Let $u_0$ be in $\mathcal{V} \cap\ell_1(w)$
and $R^N_0$ be in $\mathcal{V}^N$.
Consider $Z^N$ given by (\ref{flu}).
If $(Z^N_0)_{N\ge L}$ converges in law to $Z_0$ in $L_2(w)$ and is tight,
 then $(Z^N)_{N\ge L}$
converges in law to the unique  Ornstein-Uhlenbeck process 
solving (\ref{ousde}) starting at $Z_0$ and is tight.
\end{thm}

\section{The proof}
\setcounter{equation}{0}
\label{Sproofs} 

Let $(x)_k = x(x-1)\cdots(x-k+1)$ for $x\in \RR$ 
(the falling factorial of degree $k\in \NN$).
Let  the mappings $F^N_+$ and $F^N$ and with values in $c^0_0$
be given for $v$ in $c_0$ by 
\begin{equation}
\label{FN}
F^N_+(v) (k) = \alpha\, {(Nv(k-1))_L - (Nv(k))_L\over (N)_L}\,,
\quad
k\ge1\,,
\qquad 
F^N(v) = F^N_+(v) - F^{\vphantom{N}}_-(v)\,,
\end{equation}
where $F^{\vphantom{N}}_-$ is given in (\ref{F}).
The process $R^N$ is  Markov on $\mathcal{V}^N$, and
when in state $r$ has jumps in its $k$-th coordinate, 
$k\ge1$, of size $1/N$  at rate $N F^N_+(r) (k)$
and size $-1/N$ at rate $N F_-(r) (k)$. 

\begin{lem}
\label{dyn}
Let 
$R^N_0$ be in  $\mathcal{V}^N$,
$u$ solve (\ref{is})
starting at $u_0$ in $\mathcal{V}$,
and $Z^N$ be given by (\ref{flu}). Then
\begin{equation}
\label{zeq}
Z^N_t  = Z^N_0
+ \int_0^t \sqrt{N}\left( F^N(R^N_s) - F(u_s) \right) ds  + M^N_t
\end{equation}
defines an independent family of square-integrable martingales
$M^N = (M^N(k))_{k\in\NN}$ 
independent of $Z^N_0$ with Doob-Meyer brackets given by
\begin{equation}
\label{bracket}
\left\langle M^N(k) \right\rangle_t 
= \int_0^t \left\{ F^N_+(R^N_s)(k) + F^{\vphantom{N}}_-(R^N_s)(k) \right\} ds\,.
\end{equation}
\end{lem}

\begin{proof}
This follows from a classical application of the 
Dynkin formula.
\end{proof}

The first  lemma below shows that
it is indifferent to choose the
$L$ queues with or without replacement at this level of precision, the 
second one is a linearization formula.

\begin{lem}
\label{fact}
For $N\ge L\ge 1$ and $a$ in $\RR$ we have
\[
A^N(a) 
\defeg
{(Na)_L \over (N)_L} - a^L 
= \sum_{j=1}^{L-1} (a-1)^j a^{L-j}  
\sum_{1\le i_1<\cdots<i_j\le L-1}
{i_1\cdots i_j \over (N-i_1)\cdots (N-i_j)}
\]
and $A^N(a) =  N^{-1} O(a)$ uniformly for $a$ in $[0,1]$.
\end{lem}
\begin{proof}
We develop 
$
{(Na)_L \over (N)_L} = \prod_{i=0}^{L-1}  { Na-i \over N-i}
= \prod_{i=0}^{L-1} \left(a + (a-1){ i \over N-i}\right)
$
to obtain the identity for $A^N(a)$
and we deduce easily from it that it is
$N^{-1} O(a)$ uniformly for $a$ in $[0,1]$.
\end{proof}

\begin{lem}
\label{lin}
For $L\ge 1$ and $a$ and $h$ in $\RR$ we have
\[
B(a,h) 
\defeg
(a+h)^L - a^L - La^{L-1}h  
=\sum_{i=2}^{L} {L \choose i} a^{L-i} h^{i}
\]
with $B(a,h) =0$ for $L=1$ and $B(a,h) =h^2$
for $L=2$. For $L\ge2$ we have
$0\le B(a,h) \le h^{L} + \left(2^L -L -2\right) a h^2$
for $a$ and $a+h$ in $[0,1]$. 
\end{lem}
\begin{proof}
The identity is
Newton's binomial formula. 
A convexity argument yields $B(a,h)\ge 0$.
For
 $a$ and $a+h$ in $[0,1]$ and $L\ge2$,
$
B(a,h) \le h^L + \sum_{i=2}^{L-1} {L \choose i} a h^2
= h^L + \left(2^L - L -2\right) a h^2\,.
$
\end{proof}

For $v$ in $\mathcal{V}$ and $x$ in $c_0^0$,
considering (\ref{F}), (\ref{FN}) and Lemma~\ref{fact} let
$G^N : \mathcal{V} \rightarrow  c^0_0$ be given  by
\begin{equation}
\label{defG}
G^N(v)(k)
=\alpha A^N(v(k-1)) -  \alpha  A^N(v(k))\,,
\quad
k\ge1\,,
\end{equation}
and considering (\ref{F}), (\ref{defK}) and Lemma~\ref{lin}
let $H : \mathcal{V} \times c_0^0 \rightarrow  c_0^0$  be given by
\begin{equation}
\label{defH}
H(v,x)(k) = \alpha  B(v(k-1),x(k-1)) - \alpha B(v(k),x(k))\,,\quad 
k\ge1
\end{equation} 
so that for $v+x$ in $\mathcal{V}$ 
\begin{equation}
\label{diff}
F^N = F + G^N\,,
\qquad
F(v+x) - F(v)  
= \lin(v)  x + H(v,x)\,,
\end{equation}
and we derive the limit equation (\ref{ousde}) and (\ref{limdoob})
for the fluctuations from (\ref{zeq}) and (\ref{bracket}).

\begin{lem}
\label{l2bds}
Let $w$ satisfy (\ref{compw}).
Let $u_0$ be in $\mathcal{V} \cap\ell_1(w)$
and $R^N_0$ be in $\mathcal{V}^N$.
For $T\ge0$ we have
\[
\limsup_{N\to\infty} 
\EE\left(\left\Vert Z^N_0 \right\Vert_{L_2(w)}^2\right) < \infty
\Rightarrow
\limsup_{N\to\infty} 
\EE\biggl( \sup_{0 \le t\le T} \left\Vert Z^N_t \right\Vert_{L_2(w)}^2\biggr) 
< \infty\,.
\]
\end{lem}

\begin{proof}
Using (\ref{zeq}) and (\ref{diff})
\begin{equation}
\label{znavecg}
Z^N_t 
= Z^N_0 + M^N_t + \sqrt{N} \int_0^t G^N(R^N_s)\,ds
+\int_0^t 
\sqrt{N}\left(F (R^N_s) - F(u_s) \right)  ds
\end{equation}
where  Lemma~\ref{fact} and (\ref{compw}) yield that 
$G^N(R^N_s)(k) = N^{-1} O\!\left( R^N_s(k-1)  +  R^N_s(k) \right)$ 
and 
\begin{equation}
\label{grz}
\left\Vert G^N(R^N_s)\right\Vert_{L_2(w)}
=
N^{-1} O\!\left(\vphantom{R^N}\smash{\left\Vert R^N_s \right\Vert_{L_2(w)}}\right).
\end{equation}
We have
\begin{equation}
\label{decr1}
\left\Vert R^N_s \right\Vert_{L_2(w)}
\le
\left\Vert u_s \right\Vert_{L_2(w)} + 
N^{-1/2}  \left\Vert Z^N_s\right\Vert_{L_2(w)},
\end{equation}
Theorem~\ref{isw} yields that $F_+$, $F_-$ and $F$ are Lipschitz, the Gronwall Lemma 
that for some $K_T<\infty$
\[
\sup_{0\le t\le T}\left\Vert Z^N_t \right\Vert_{L_2(w)}
\le
K_T\biggl(\left\Vert Z^N_0 \right\Vert_{L_2(w)} 
+ N^{-1/2} K_T \left\Vert u_0 \right\Vert_{L_2(w)}
+ \sup_{0\le t\le T}\left\Vert M^N_t \right\Vert_{L_2(w)}
\biggr),
\]
and we conclude using the Doob inequality, (\ref{bracket}), (\ref{diff}),
\begin{equation}
\label{brkbd}
\left\Vert F_+(R^N_s) + F_-(R^N_s) \right\Vert_{L_2(w)}
\le
K\left\Vert R^N_s \right\Vert_{L_2(w)},
\end{equation}
and the bounds (\ref{grz}) and (\ref{decr1}).
\end{proof}

\begin{lem}
\label{tight}
Let $w$ satisfy (\ref{compw}), and consider 
$L_2(w)$
with its weak topology and
$\DD(\RR_+, L_2(w))$ with the corresponding Skorokhod
topology.
Let $u_0$ be in $\mathcal{V} \cap\ell_1(w)$
and $R^N_0$ be in $\mathcal{V}^N$.
Consider $Z^N$ given by (\ref{flu}).
If $(Z^N_0)_{N\ge L}$ is tight
 then $(Z^N)_{N\ge L}$ is tight and its limit points are continuous.
\end{lem}

\begin{proof}
For $\varepsilon>0$ let $r_\varepsilon < \infty$ be such that
$\PP(Z^N_0 \in B(r_\varepsilon))> 1 - \varepsilon$ for $N\ge1$
(see the discussion prior to Theorem~\ref{fclt}).
Let $R^{N,\varepsilon}_0$ be equal to $R^N_0$ 
on $\{Z^N_0 \in B(r_\varepsilon)\}$ and such that
$Z^{N,\varepsilon}_0$ is uniformly
bounded in $L_2(w)$ on
$\{Z^N_0 \not\in B(r_\varepsilon)\}$ (for instance deterministically equal to
some outcome of
$R^N_0$ on $\{Z^N_0 \in B(r_\varepsilon)\}$).
Then $Z^{N,\varepsilon}_0$ is uniformly bounded in $L_2(w)$ and we may
 use a coupling argument
to construct $Z^{N,\varepsilon}$ and $Z^N$ coinciding
on $\{Z^N_0 \in B(r_\varepsilon)\}$.

Hence to prove tightness of $(Z^N)_{N\ge L}$ we may restrict our attention
to $(Z^N_0)_{N\ge L}$ uniformly bounded in $L_2(w)$, for which we may use 
Lemma \ref{l2bds}.

The compact subsets of $L_2(w)$ are Polish, a fact yielding tightness criteria. 
We deduce from Theorems~4.6 and 3.1 in 
Jakubowski~\cite{Jakubowski:86}, which considers
completely regular Hausdorff spaces (Tychonoff spaces)
of which $L_2(w)$ with its weak topology is an example,
that $(Z^N)_{N\ge L}$ is tight if
\begin{enumerate}
\item
For each $T\ge 0$ and $\varepsilon >0$ there is a bounded subset
$K_{T,\varepsilon}$ of $L_2(w)$ such that for
$N\ge L$ we have
$
\PP\!\left( Z^{N} \in \DD([0,T], K_{T,\varepsilon}) \right) > 1-\varepsilon
$.

\item
For each $d\ge1$, the 
$d$-dimensional processes $(Z^{N}(1), \ldots, Z^{N}(d))_{N \ge L}$
are tight.
\end{enumerate}

Lemma~\ref{l2bds} and the  Markov inequality yield condition 1.
We use (\ref{znavecg}) (derived from (\ref{zeq}))
 and (\ref{bracket}) and (\ref{diff}), and the  bounds (\ref{grz}), 
(\ref{decr1}) and (\ref{brkbd}). The uniform bounds in Lemma~\ref{l2bds}
and the fact that $Z^{N}(k)$ has jumps of size $N^{-1/2}$ classically imply
 that the above finite-dimensional processes are tight 
and have continuous limit points, 
see for instance  Ethier-Kurtz~\cite{Ethier:86} Theorem~4.1 p.~354
or Joffe-M{\'e}tivier~\cite{Joffe:86} Proposition~3.2.3 
and their proofs.
\end{proof}

{\em End of the proof of Theorem \ref{fclt}\/}.
Lemma \ref{tight} implies that from any subsequence of $Z^N$
we may extract a further subsequence which converges to
some $Z^\infty$ with continuous sample paths. Necessarily
$Z^\infty_0$ has same law as $Z_0$. In  (\ref{znavecg}) we have
considering (\ref{diff}) 
\begin{equation}
\label{fin}
\sqrt{N}\!\left(F (R^N_s)(k) - F(u_s)(k) \right)
= \lin(u_s) Z^N_s + \sqrt{N} H\!\left(u_s,\smash{N^{-1/2}}Z^N_s\right).
\end{equation}
We use the  bounds (\ref{grz}), (\ref{decr1}) and (\ref{brkbd}),
the uniform bounds in Lemma~\ref{l2bds},
and additionally (\ref{defH}) and Lemma~\ref{lin}. 
We deduce by a martingale characterization
that $Z^\infty$ has the law of the Ornstein-Uhlenbeck process
unique solution for (\ref{ousde})
in $L_2(w)$
starting at $Z^\infty_0$, see Theorem~\ref{Kbdd}; 
the drift vector is given by the limit for (\ref{zeq}) and
(\ref{znavecg}) considering (\ref{fin}),
and the martingale bracket by the limit for (\ref{bracket}).
See for instance 
Ethier-Kurtz~\cite{Ethier:86} Theorem~4.1 p.~354 
or Joffe-M{\'e}tivier~\cite{Joffe:86} Theorem 3.3.1
and their proofs for details. Thus, this law is the unique
accumulation point for the relatively compact sequence of laws of $(Z^N)_{N\ge1}$,
which must then converge to it, proving Theorem \ref{fclt}.


\end{document}